\documentclass{article}
\usepackage[latin1]{inputenc}
\usepackage{amsmath, amsthm, amssymb}

\newcommand{\ra}{\mathop{\rightarrow }}

\newcommand{\cA}{{\cal A}}
\newcommand{\cB}{{\cal B}}
\newcommand{\cC}{{\cal C}}

\newcommand{\cF}{{\cal F}}

\newcommand{\cN}{{\cal N}}
\newcommand{\cO}{{\cal O}}

\newcommand{\bbZ}{\mathbb Z}
\newcommand{\bbR}{\mathbb R}
\newcommand{\bbN}{\mathbb N}

\newcommand{\bbQ}{\mathbb Q}
\newcommand{\bbE}{\mathbb E}
\newcommand{\bbT}{\mathbb T}
\newcommand{\bbP}{\mathbb P}

\newtheorem{theo}{Theorem}[section]
\newtheorem{pr}{Proposition}[section]

\newtheorem{lem}{Lemma}[section]
\newtheorem{defn}{Definition}[section]

\begin{document}
\hoffset=-1,5cm
\parskip=4mm
\title{The Curie-Weiss model with dynamical external field}
\author{C. Dombry and N. Guillotin-Plantard \thanks{Universit\'e Claude Bernard-Lyon 1, Institut Camille Jordan,
43 bld du 11 novembre 1918, 69622 Villeurbanne, France, e-mails: dombry@univ-lyon1.fr, nadine.guillotin@univ-lyon1.fr}\protect\hspace{1cm}}
\date{}
\maketitle
~\\
Key words: Curie-Weiss model, dynamic random walk, dynamical system, ergodic theory, diophantine approximations, 
limit theorems, large deviation principle, statistical physics.\\
~\\
AMS Subject classification: Primary: 60K35; secondary: 82B41; 82B44; 60G50

\begin{abstract}
We study a Curie-Weiss model with a random external field generated by a dynamical system. Probabilistic limit theorems (weak law of large numbers, 
central limit theorems) are proven for the corresponding magnetization. Our results extend those already obtained in \cite{livEll} and \cite{ell1}.
\end{abstract}

\section{Introduction}
The Curie-Weiss model is a well-known approximation to the Ising model (see \cite{livEll}). Probabilistic limit theorems for the Curie-Weiss model
have been proven by the following authors: Ellis and Newman \cite{ell1}, Ellis, Newman and Rosen \cite{ell2}...
The purpose of the present paper is to prove limit theorems for the Curie-Weiss model with random external field generated by a dynamical system, namely weak law of large numbers and central limit theorems for the associated magnetization. Our main motivation is the understanding of the statistical properties of the following physical model: consider $n$ particles (iron atoms for instance) distributed along a lattice $\Gamma=\{1,\ldots, n\}$. The value of $\pm 1$ at a site represents the {\it spin}, or {\it magnetic moment}, of the particle at that site. The particles are placed in a magnetic field which is given in terms of a dynamical system $S=(E,{\cal A},\mu,T)$, that is a probability space $(E,{\cal A},\mu)$, $T$ a transformation of $E$ and a function $f$ defined on $E$ with values in $[0,1]$. Let $\beta >0$ be the inverse temperature and $J$ a coupling constant assumed strictly positive.\\*
Given a configuration  $\sigma=(\sigma_i)_{i=1\ldots,n} \in \Omega_n=\{-1,+1\}^{n}$ and $x\in E$, we define 
the Hamiltonian,
$$ H_{n,x}(\sigma)=  \frac{\beta J}{2n} \left(\sum_{i=1}^n \sigma_i\right)^2  +\frac{1}{2} \sum_{i=1}^n \log\left(\frac{f(T^i x)}{1-f(T^i x)}\right) \sigma_i.$$
We denote by $\bbQ_{n,x}$ the Gibbs measure on  $\Omega_n$ defined by
$$ \bbQ_{n,x}(\sigma) = \frac{1}{Z_{n,x}}\exp[H_{n,x}(\sigma)]$$
where $Z_{n,x}$ is the normalizing constant, called partition function
$$Z_{n,x}= \sum_{\sigma\in \Omega_n} \exp[H_{n,x}(\sigma)].$$
Equivalently,  $\bbQ_{n,x}$ is the probability measure defined on $(\Omega_n, {\cal P}(\Omega_n))$ such that, for any $A\in {\cal P}(\Omega_n) $, 
$$ \bbQ_{n,x}(A) = \frac{1}{\tilde{Z}_{n,x}}\int_{A} \exp[H_{n,x}(\sigma)] \, {\rm d}\bbP_n(\sigma)$$
where  $\bbP_n$ is the uniform distribution on $\Omega_n$ (i.e. for any $\sigma\in \Omega_n$, $\bbP_n(\sigma) =\frac{1}{2^n}$)
and $\tilde{Z}_{n,x}$ is the partition function 
$$\tilde{Z}_{n,x}= \int_{\Omega_n} \exp[H_{n,x}(\sigma)] \, {\rm d}\bbP_n(\sigma).$$
For each configuration $\sigma=(\sigma_i)_{i=1,\ldots,n}$ we define the associated {\it magnetization} (or {\it total spin})
$$M_n=\sum_{i=1}^{n} \sigma_i.$$
Remark that when $f\equiv 1/2$ our model corresponds to the Curie-Weiss model without external field studied in \cite{livEll, ell1, ell2}.
Moreover it is worth remarking that any field $(g(T^i x))_{i\geq 1}$ can be considered by choosing the function $f$ as $e^{g}/(1+e^{g})$. In particular, it includes the case where the field is given in terms of a sequence of independent and identically Bernoulli random variables taking the values $-\varepsilon$ and $\varepsilon$ with probability $1/2$ considered in \cite{Kul} and \cite{Picco}. We will be mostly interested in the special case when the dynamical system is the irrational rotation on the torus which corresponds to a quasiperiodic random field; we refer to \cite{KPZ} for a complete and precise discussion about the relevance of this model in the modelization of certain physical models.\\*
We are interested in studying the asymptotic behaviour of
$M_n$ in the so-called thermodynamical limit $n\rightarrow +\infty$. In \cite{ell2} a physical interpretation of this limit behaviour is given in relation with stable states (mixed or pure) and metastable states of the underlying physical system. An illustrative example derived from thermodynamics is given, namely a detailed description of states as well as the phase transition in a gas-liquid system.\\*
At infinite temperature (i.e. $\beta=0$) the probability measure $\bbQ_{n,x}$ is equal to the product measure 
$$ \prod_{i=1}^n \left(f(T^ix)\delta_1+ (1-f(T^ix))\delta_{-1}\right)$$
This implies that the random variable $M_n$ is just a sum of independent random variables $\sigma_i$ taking the value $1$ with probability $f(T^{i}x)$ and $-1$ with probability $1-f(T^{i}x)$. So, in this particular case, the sequence of random variables $(M_n)_{n\geq 1}$ is a so-called {\it dynamic $\bbZ-$random walk} (see Section 2). The dynamic $\bbZ$-random walks were introduced by the second author in \cite{gui1}, then generalized to dimension $d>1$ in \cite{guiJTP}. 
Theoretical results about dynamic random walks and their applications can be found in the recent book \cite{guiliv}. We are mainly interested in limit theorems (i.e. strong law of large numbers, central limit theorem and large deviation principle) for dynamic $\bbZ$-random walks. We recall some of them here under simplified assumptions: assume that $E$ is a compact metric space, ${\cA}$ the associated Borel $\sigma$-field, $T$ a continuous transformation of $E$. If there exists an unique invariant measure
$\mu$ i.e. $(E,{\cA},\mu,T)$ is uniquely ergodic and if $f$ is continuous with integral equal to $1/2$, then for every $x\in E$, $(M_n/n)_{n\geq 1}$ converges almost surely to 0 as $n$ goes to infinity. Moreover, if we assume that $a=\int_E 4f(1-f)\, d\mu>0$ and that
\begin{equation}\label{h}
\sum_{i=1}^n [f(T^i x) - 1/2]= o(\sqrt{n})
\end{equation}
then the sequence $(M_n/\sqrt{n})_{n\geq 1}$ converges in distribution as $n\rightarrow +\infty$ to the Normal distribution ${\cal N}(0,a)$. Let us recall that an important feature of certain Ising model is the existence of a critical value $\beta_c$ of $\beta$:   For $0<\beta<\beta_c$, the spins are weakly correlated and the probabilistic limit theorems obtained at $\beta=0$ are valid. For $\beta>\beta_c$, the correlation between the spins is strongly positive and the limit results are completely different. The model is then said to present a phase transition at $\beta=\beta_c$.       
We will prove for our model that the limit theorems obtained at $\beta=0$ for the dynamic random walk are still valid for any $\beta<\beta_c$ under the same hypotheses. The critical value $\beta_c$ is shown to fluctuate between $1/J$ and $1/(Ja)$ according to the dynamical system and the function $f$ we consider. For an explicit class of dynamical systems and functions $f$, we are able to prove that $\beta_c$ is equal to $1/(Ja)$ and that at $\beta=\beta_c$, under suitable assumptions on $f$, there exists some $\gamma\in (0,1)$ so that as $n\rightarrow +\infty$, $(M_n/n^{\gamma})_n$ converges in distribution to an explicit non Gaussian random variable. For $\beta>\beta_c$, the situation is not so well understood. Let us recall that for the Curie-Weiss model with zero external field (e.g. $f\equiv 1/2$), the sequence $(M_n/n)_n$ converges in distribution to
$1/2(\delta_m+\delta_{-m})$ where $m>0$ (see for instance Theorem IV.4.1 in \cite{livEll}). 
In Section 4 we show that if the function $f$ is not identically equal to $1/2$, when $\beta>\beta_c$, the sequence $(M_n/n)_n$ does not converge in distribution. We conjecture that Theorem 1 in \cite{Kul} should be true for our model, under suitable assumptions; this work is in progress. \\*
In \cite{KPZ}, the authors extended the Pirogov-Sinai theory to a class of models with small quasiperiodic interactions as perturbations of the periodic ones. More precisely, the low temperature phase diagram for spin systems with periodic hamiltonians perturbated by quasiperiodic interactions is studied. Under diophantine conditions and derivability conditions on the interaction potentials they prove that the low temperature phase diagram is a homeomorphic deformation of the phase diagram at zero temperature. Our model in the case when the dynamical system $S$ is an irrational rotation on the torus belongs to this class of hamiltonians perturbated by quasiperiodic ones. In Section 5 the same kind of conditions on the diophantine approximation of the irrational angle and on the smoothness of the function $f$ will be needed in order to state the limit theorems.   
\\*   

The outline of the paper is as follows: In Section 2, we define the dynamic $\bbZ$-random walk and recall some results which will be useful in the sequel.
In Section 3, we state and prove our results under general assumptions. In Section 4, we apply results of Section 3 when the integral of $f$ is equal to $1/2$ and $a>0$. In Section 5, we study the case when the dynamical system is given by an irrational rotation on the torus.

\section{Dynamic $\bbZ$-random walks}
The dynamic random walks were introduced in \cite{gui1} and generalized to upper dimensions in \cite{guiJTP}. We now recall
some of the results obtained in dimension one. 
Let $S=(E,\cA,\mu,T)$ be a dynamical system where $(E,\cA,\mu)$ is
a probability space and $T$ is a measure-preserving transformation
defined on $E$. Let
$f$ be a measurable function defined on $E$ with values in $[0,1]$. For each $x\in E$, we denote by $\bbP_{x}$  the distribution of the time-inhomogeneous random walk:
$$S_{0}=0,\ \ S_{n}=\sum_{i=1}^{n} X_{i}~\ \mbox{for}~~ n\geq 1$$
with step distribution 
\begin{equation}\label{step}
\bbP_{x}(X_{i}=z)=\left\{\begin{array}{ll}
 f(T^{i}x) & \mbox{if}\ \ z=1\\
1-f(T^{i}x) & \mbox{if}\ \ z=-1\\
0 & \mbox{otherwise.}
\end{array} \right.
\end{equation}
The expectation with respect to $\bbP_{x}$ will be denoted by $\bbE_x$.
It is worth remarking that if the function
$f$ is not constant, $(S_{n})_{n\in \bbN}$ is a
non-homogeneous Markov chain. This Markov chain can be classified
in the large class of random walks evolving in a random
environment. In most of the papers (see for instance \cite{GDH},
\cite{com},...), the environment field takes place in space but it
can also take place in space and time (see \cite{pel}). Following
the formalism used in the study of these random walks, when $x$ is
fixed, the measure $\bbP_{x}$ is called \mbox{\it quenched} and
the measure averaged on values of $x$ defined as $\bbP(.)=
\int_{E}\bbP_{x}(.)\ {\rm d}\mu(x)$ is called \mbox{\it annealed}. 
In the quenched case, the random variables $X_i, i\geq 1$ are independent, but not necessarily identically distributed. In the annealed case, the $X'$s defines a stationary sequence of dependent random variables whose the correlations are related to the ones of the underlying dynamical system. We refer to \cite{arno} (Section 2.1) for a more precise discussion of these two cases.\\*
Let ${\cC}_{1/2}(S)$ denote the class of functions $f\in L^{1}(\mu)$ satisfying the following condition: for every $x\in E$,
$$\left|\sum_{i=1}^{n}\Big(f(T^{i}x)-\int_{E}f d\mu\Big)\right|=o\Big(\sqrt{n}\Big).$$
Let us assume that $f\in{\cC}_{1/2}(S)$; if $\int_{E}f\,  {\rm d}\mu=\frac{1}{2}$ and $a=\int_E 4f(1-f)\, {\rm d}\mu>0$, then, for every $x\in E$, the sequence of random variables $(S_n/\sqrt{n})_{n\geq 1}$
 converges in distribution to the Normal law ${\cal N}(0,a)$ (see \cite{GPS}).  A
strong law of large numbers for the dynamic $\bbZ$-random walk can be
obtained for $\mu$-almost every point $x\in E$ from Kolmogorov's
theorem assuming that the function $f$ is measurable (see Chapter 2 in \cite{guiliv} for details). The limit is then given by $2\, \bbE(f|{\cal
I})-1$ where ${\cal I}$ is the invariant
$\sigma$-field associated to the transformation $T$. So,
$(S_{n}/n)_{n\ge 1}$ is a good candidate for a large deviation principle. Let us recall what is a large deviation principle:
Let $\Gamma$ be a Polish space endowed with the Borel
$\sigma$-field ${\cal B}(\Gamma)$. A good {\it rate} function is a
lower semi-continuous function $\Lambda^{*}: \Gamma \ra
[0,\infty]$ with compact level sets $\{x; \Lambda^{*}(x)\leq
\alpha\}, \alpha\in [0,\infty[.$ Let $v=(v_{n})_{n}\uparrow\infty$
be an increasing sequence of positive reals. A sequence of random
variables $(Y_{n})_{n}$ with values in $\Gamma$ defined
 on a probability space $(\Omega, {\cF}, \bbP)$ is said to satisfy
 {\it a Large Deviation Principle} (LDP)
 with speed $v=(v_{n})_{n}$ and the good rate function $\Lambda^{*}$ if
for every Borel set $B\in{\cal B}(\Gamma)$,
\begin{eqnarray*}
-\inf_{x\in B^{o}}\Lambda^{*}(x)&\leq&\liminf_{n}\frac{1}{v_{n}}\log \bbP(Y_{n}\in B)\\
&\leq&\limsup_{n}\frac{1}{v_{n}}\log \bbP(Y_{n}\in B)\leq -\inf_{x\in\bar{B}}\Lambda^{*}(x).
\end{eqnarray*}
\begin{theo}\label{th1}
\begin{enumerate}
\item  For $\mu$-almost every $x\in E$, the sequence $(S_{n}/n)_n$ satisfies a LDP with
speed $n$ and good rate
  function
  $$\Lambda^{*}(y)=\sup_{\lambda\in\bbR} \{<\lambda,y>-\Lambda(\lambda)\}$$
  where
  $$\Lambda(\lambda)=\bbE\left(\log\Big(e^{\lambda} f +
  (1-f) e^{-\lambda}\Big)\Big|\ {\cal I}\right),$$
${\cal I}$ being the $\sigma$-field generated by the fixed points of the transformation $T$.
\item Assume that $E$ is a compact metric space, ${\cA}$ the associated Borel $\sigma$-field
and $T$ a continuous transformation of $E$. If $(E,{\cA},\mu,T)$ is uniquely ergodic (i.e. there exists an unique invariant measure
$\mu$) and if $f$ is continuous, then 1. holds for every $x\in E$.\\*
The rate function is then deterministic and equal to
$$\Lambda(\lambda)=\int_{E} \log\Big(e^{\lambda} f(x) + (1-f(x)) e^{-\lambda} \Big)\, {\rm d}\mu(x) .$$
\end{enumerate}
\end{theo}
\noindent Let us mention that an annealed large deviations
statement for $(S_n/n)_n$ under the measure $\bbP$ can easily be proved using results of \cite{din} (Remark that $E$ is assumed to be compact).
The proof of the above theorem can be found in \cite{DGPSP}.

\section{Limit theorems for the magnetization}
In this section, $E$ is assumed to be a compact metric space, ${\cA}$ the associated Borel $\sigma$-field, $\mu$ a probability measure on $(E,{\cA})$, $T$ a continuous measure-preserving transformation of $E$, $f$ a continuous function from $E$ to $[0,1]$ and $x$ a fixed point of $E$. The system $(E,{\cA},\mu,T)$ is asumed to be uniquely ergodic. In the sequel, the sequence $(S_n)_n$ will denote the dynamical random walk introduced in Section 2 and $(M_n)_n$ the magnetization defined in the introduction. 

\subsection{Weak law of large numbers for the magnetization}
For every $n\geq 1$, we define the function 
\begin{eqnarray*}
G_n(s) &=&\frac{\beta J}{2} s^2 - \frac{1}{n} \log\ \bbE_x(\exp(\beta J s S_n))\\
 &=& \frac{\beta J}{2} s^2 - \frac{1}{n}\sum_{i=1}^n  L(f(T^i x),\beta J s) 
\end{eqnarray*}
where the function $L$ is defined on $[0,1]\times\bbR$ by
$$L(\phi,s)= \log\left( \phi\ e^{ s}  + (1-\phi)\ e^{- s}\right ) .$$
We also define the function $G$ by 
$$G(s)= \frac{\beta J}{2} s^2 - \int_E  L(f(y),\beta Js)  \, {\rm d}\mu(y).$$

\begin{theo}\label{theo1}
The function $G$ is real analytic, and the set where $G$ achieves its minimum is non-empty and finite.
\end{theo}

\begin{defn} We will denote by $g=\min\{G(s);s\in\bbR\}$ the minimum of $G$ and by $m_1,\cdots,m_r$ the points where $G$ is minimal. Let us define the type $2k_i$  and the strength $\lambda_i>0$ of the minimum $m_i$ by
\begin{eqnarray*}
2k_i&=&\min\{j\geq 1\ |\  G^{(j)}(m_i)\neq 0 \},\\
\lambda_i&=&G^{(2k_i)}(m_i).
\end{eqnarray*}
\end{defn}
\noindent Let us remark that $g$ is nonpositive since $G(0)=0$.

For every $\alpha\in [0,1]$, we define $\cC_{\alpha}(S)$ the class of $\mu$-integrable functions 
$h:E\rightarrow \bbR $ satisfying the following condition: for every point $x\in E$,
$$\left|\sum_{k=1}^{n}\Big(h(T^{k}x)-\int_{E}h \,{\rm d}\mu\Big)\right|=o\Big(n^{\alpha} \Big).$$
Remark that since the dynamical system is uniquely ergodic, the class $\cC_{1}(S)$ always contains the set of continuous functions on $E$.
\begin{theo}\label{theo2}
\begin{enumerate}
\item
Assume that for every $i\in\{1,\ldots,r\}$ and every $j\in\{1,\ldots,2k_i\}$, 
the function 
$$\begin{array}{lll}
E&\rightarrow& \bbR \\
y&\mapsto& \frac{\partial ^j}{\partial s^j}L(f(y),\beta J m_i)
\end{array}$$ 
belongs to the set $\cC_{j/2k_i}(S).$\\
Then, for every bounded continuous function $h$, the expectation  of $h(M_n/n)$ under $\bbQ_{n,x}$ is equivalent, as $n$ goes to infinity, to
$$\displaystyle\frac{\displaystyle\sum_{i=1}^r b_{i,n} h(m_i)}{\displaystyle\sum_{i=1}^r b_{i,n}}$$
where 
$$b_{i,n}=n^{-1/2k_i}e^{-nG_n(m_i)}\lambda_i^{-1/2k_i}\int_{-\infty}^{+\infty} \exp(-s^{2k_i}/(2k_i)!)\, {\rm d}s.$$
In particular, if $G$ achieves its minimum at a unique point $m$, then the distribution of $M_n/n$ under $\bbQ_{n,x}$ converges to $\delta_m$ the Dirac mass at $m$.
\item The distribution of $M_n/n$ under $\bbQ_{n,x}$ verifies a large deviation principle with speed $n$
  and good rate function
  $$I_{\beta,x}(z)=  \Lambda^{*}(z) - \frac{\beta J}{2}z^2 - \inf_{z\in\bbR} \{ \Lambda^{*}(z) - \frac{\beta J}{2} z^2\}$$
  where $\Lambda^{*}$ is defined in Theorem \ref{th1}.
\end{enumerate}
\end {theo}

\subsection{Scaling limit for the magnetization}

\begin{theo}\label{theo3}
Assume that $G$ has a unique global minimum $m$ of type $2k$ and strength $\lambda$ 
and that for every $j\in\{1,\ldots,2k\}$, 
the function $\frac{\partial ^j}{\partial s^j}L(f(.),\beta J m)$ belongs to the set $\mathcal{C}_{j/2k}(S).$
Then, the following convergence of measures holds:
$$\frac{M_n-nm}{n^{1-1/2k}}\Rightarrow Z(2k,\tilde \lambda) $$
where $Z(2k,\tilde \lambda)$ is the probability measure with density function
$$C \exp\left(-\tilde\lambda s^{2k}/(2k)!\right),$$
where $C$ is a normalizing constant and $\tilde\lambda$ is defined by
$$\tilde\lambda =\left\{\begin{array}{lll}
\left(\frac{1}{\lambda}-\frac{1}{\beta J}\right)^{-1} &\ {\rm\ if}\ & k=1\\ 
\ \ \ \ \ \lambda &\ {\rm\ if}\ & k\geq 2\end{array}\right..$$
\end{theo} 
\smallskip
\noindent{\bf Remark:}
Note that the case of a minimum of type 2 yields a central limit theorem: the fluctuations of $M_n/n$ around $m$ are of order $n^{-1/2}$ and Gaussian.
When the type of the minimum is greater than 4, the limit distributions are non standard.
\subsection{Technical lemmas}
\begin{lem} \label{lemma1}
Let $Y$ be a random variable with distribution ${\cal N}(0,1/(\beta J))$, independent of $M_n$ for every $n\geq 1$. Then, given $m$ and $\gamma$ real, the probability density function of the random variable 
$$\frac{Y}{n^{1/2-\gamma}}+\frac{M_n-nm}{n^{1-\gamma}} $$ 
is equal to
$$\frac{ \exp{(-nG_n(m+sn^{-\gamma}))} }{\int_{\bbR}\exp{(-nG_n(m+sn^{-\gamma})) }\, {\rm d}s}.$$
\end{lem}
\noindent{\bf Proof:}\\*
The probability density function of the random variable ${n}^{1/2} Y+ M_n $ is given by
$$ \frac{1}{Z_n} \sqrt{\frac{\beta J}{2\pi n}} \int e^{-\beta J(s-x)^2/(2n)} e^{\beta Jx^2/2n}\,  {\rm d}\tilde{\rho}_n(x)$$
with $\tilde{\rho}_n= \rho_1*\ldots *\rho_n$ where $\rho_i=f(T^i x) \delta_1+(1-f(T^i x)) \delta_{-1}$. It can be rewritten as
$$ \frac{1}{Z_n} \sqrt{\frac{\beta J}{2\pi n}} e^{-\beta Js^2/2n} \bbE_x(e^{\beta Js S_n/n})$$
So, by a change of variables, the probability density function of the random variable $\frac{Y}{n^{1/2-\gamma}}+\frac{M_n-nm}{n^{1-\gamma}} $  is given by
$$ \frac{1}{Z_n} \sqrt{\frac{\beta J}{2\pi }}n^{1/2-\gamma} e^{-\beta Jn(m+sn^{-\gamma})^2/2} \bbE_x(e^{\beta J (m+sn^{-\gamma}) S_n})$$
The lemma is then easily deduced.

The previous lemma suggest that the behaviour of the sequence of random variables $M_n$ and of the sequence of functions $G_n$ are linked together. 
\begin{lem}\label{lemma2}
The sequence of functions $(G_n)_{n\geq 1}$ converges to $G$ uniformly on compacta of $\bbR$ as $n$ goes to infinity. 
Furthermore, for every $k\geq 1$, the sequence of derivative functions $(G_n^{(k)})_{n\geq 1}$ converges to $G^{(k)}$ uniformly on compacta of $\bbR$ as $n$ goes to infinity.
\end{lem}
\noindent{\bf Proof:}\\*
Note that the function $L$ is of class $\mathcal{C}^{\infty}$ on $[0,1]\times\bbR$, so that for any $s\in\bbR$, the function $y \mapsto \frac{\partial^k}{\partial s^k} L(f(y),\beta J s)$ is continuous.
For any $s\in\bbR$ and $k\geq 0$,  
$$G_n^{(k)}(s)-G^{(k)}(s)=(\beta J)^k \left[\frac{1}{n}\sum_{i=1}^n \frac{\partial^k}{\partial s^k} L(f(T^ix),\beta J s)-\int_E \frac{\partial^k}{\partial s^k}L(f(y),\beta J s)  {\rm d} \mu(y)\right]$$
The unique ergodicity hypothesis implies that this quantity converges to $0$ as $n$ goes to infinity.\\  
We prove the uniform convergence with the following majoration of the difference $|G_n^{(k)}(s)-G^{(k)}(s)|$ 
on the compact $[-\alpha;\alpha]$: 
$$|G_n^{(k)}(s)-G^{(k)}(s)|\leq |G_n^{(k)}(0)-G^{(k)}(0)|+\int_{-\alpha}^\alpha |G_n^{(k+1)}(s)-G^{(k+1)}(s)| \, {\rm d}s.$$
The unique ergodicity hypothesis implies that for any $s\in\bbR$, $|G_n^{(k+1)}(s)-G^{(k+1)}(s)|$ converges to zero as $n$ goes to infinity. Furthermore, the function $\frac{\partial^k}{\partial s^k}L$ is bounded on the compact set $[0,1]\times[-\beta J\alpha,\beta J\alpha]$, and hence the difference $|G_n^{(k+1)}(s)-G^{(k+1)}(s)|$ is uniformly bounded on $[-\alpha,\alpha]$.
 Finally, by dominated convergence theorem,
$$|G_n^{(k)}(s)-G^{(k)}(s)|\leq |G_n^{(k)}(0)-G^{(k)}(0)|+\int_{-\alpha}^\alpha |G_n^{(k+1)}(s)-G^{(k+1)}(s)| \, {\rm d}s \rightarrow 0,$$
and the convergence is uniform for $s\in[-\alpha,\alpha]$.

\begin{lem}\label{lemma3}
Let $m$ be a global minimum of $G$ of type $2k$ and strength $\lambda$.
Suppose that for every $j\in\{1,\ldots,2k\}$, the function $y\mapsto\frac{\partial^j }{\partial s^j}L(f(y),\beta J m)$ belongs to the class $\cC_{\frac{j}{2k}}(S)$. Then, for every $s\in \bbR$,
\begin{equation}\label{DL1}
\lim_{n\rightarrow +\infty} n \Big(G_n\left(m+sn^{-1/2k}\right)-G_n(m)\Big)=\lambda\frac{s^{2k}}{(2k)!}.
\end{equation}
Furthermore, there exist $\delta>0$ and $N\geq 1$  such that for every $n\geq N$ and $s\in [-\delta n^{1/2k};\delta n^{1/2k}]$,
\begin{equation}\label{DL2}
n \left(G_n\left(m+sn^{-1/2k}\right)-G_n(m)\right)\geq  \frac{\lambda}{2}\frac{s^{2k}}{(2k)!}-\sum_{j=1}^{2k-1}|s|^j.
\end{equation}
\end{lem}
\noindent{\bf Proof:}\\*
Let $s\in\bbR$ and $u=sn^{-1/2k}$. Taylor's formula implies that 
$$G_n(m+u)-G_n(m)=\sum_{j=1}^{2k} \frac{G_n^{(j)}(m)}{j!}u^j+R_n(u),$$
where the remainder $R_n$ has the integral form
$$R_n(u)=\frac{u^{2k+1}}{(2k)!}\int_0^1 (1-\theta )^{2k}G_n^{(2k+1)}(m+\theta u)\ {\rm d}\theta.$$
The $j$-th derivative of $G_n$ at point $m$ is equal to
$$G_n^{(j)}(m)=P_j(m)-\frac{(\beta J)^j}{n}\sum_{i=1}^n \frac{\partial^j }{\partial s^j}L(f(T^i x),\beta Jm),$$
where $$P_j(m)=\left\{\begin{array}{ll} \beta J m & {\rm \ if\ } j=1 \\ \beta J &{\rm \ if\ } j=2 \\ 0& {\rm \ otherwise} \end{array}\right..$$
As $n$ goes to infinity, this quantity converges to 
$$G^{(j)}(m)=P_j(m)-(\beta J)^j\int_E \frac{\partial^j }{\partial s^j}L(f(y),\beta J m) \,{\rm d}\mu(y).$$
The hypothesis that the function $y\mapsto\frac{\partial^j }{\partial s^j}L(f(y),\beta J m)$ belongs to the class $\cC_{\frac{j}{2k}}(S)$ implies that for every $j\in\{1,\ldots,2k\}$,
$$n|G_n^{(j)}(m)-G^{(j)}(m)|=o\Big(n^{\frac{j}{2k}}\Big).$$
Since the point $m$ is a global minimum of $G$ of type $2k$ and strength $\lambda$,  $G^{(j)}(m)=0$ for every $j\in\{1,\ldots,2k-1\}$ and $G^{(2k)}(m)=\lambda>0$. This implies that as $n$ goes to infinity, for every $j\in\{1,\ldots,2k-1\}$,
\begin{equation}\label{ASY1}
G_n^{(j)}(m)n^{1-\frac{j}{2k}}\rightarrow 0,
\end{equation}
and that for $j=2k$,
\begin{equation}\label{ASY2}
G_n^{(2k)}(m)\rightarrow \lambda.
\end{equation}
The integral remainder satisfies 
$$nR_n(sn^{-1/2k})=\frac{s^{2k+1}n^{-1/2k}}{(2k)!}\int_0^1 (1-\theta )^{2k}G_n^{(2k+1)}(m+\theta sn^{-1/2k})\ {\rm d}\theta .$$
The fact that the functions  $G_n^{(2k+1)}$ are uniformly bounded on the compact set $[m-|s|, m+|s|]$ implies that the integral is bounded as $n$ goes to infinity, so that
$$nR_n(sn^{-1/2k})\rightarrow 0.$$
Hence the formula 
$$n \Big(G_n\left(m+sn^{-1/2k}\right)-G_n(m)\Big)=
\sum_{j=1}^{2k} \frac{G_n^{(j)}(m)}{j!} n^{1-j/2k} s^j + nR_n(sn^{-1/2k}),$$
yields the limit $\lambda s^{2k}/(2k)!$ as $n$ goes to infinity and this proves equation (\ref{DL1}).\\
Let us prove (\ref{DL2}), from (\ref{ASY1}) and (\ref{ASY2}), there exists $N\geq 1$ such that for every $n\geq N$, and every $j\in \{1,\ldots,2k-1\}$, 
$$\left| \frac{G_n^{(j)}(m)n^{1-j/2k}}{j!} \right|\leq 1$$
and, for $j=2k$,
$$G_n^{(2k)}(m)\geq 3\lambda /4.$$
There also exists  $\delta>0$ such that for every $n\geq 1$, and every $s\in[-\delta n^{1/2k};\delta n^{1/2k}]$, 
$$\Big| sn^{-1/2k}\int_0^1 (1-\theta )^{2k}G_n^{(2k+1)}(m+\theta sn^{-1/2k})\ {\rm d}\theta \Big| \leq \lambda /4,$$
which implies that for every $n\geq 1$ and every $s\in [-\delta n^{1/2k};\delta n^{1/2k}]$,
$$nR_n(sn^{-1/2k})\geq -\frac{\lambda}{4} \frac{s^{2k}}{(2k)!}.$$
(There exists $M$ such that $|G_n^{(2k+1)}(s)|\leq M$ for every $n\geq 1$ and every $s\in [m-1,m+1]$, then choose $\delta$ as the minimum of $(2k+1)\lambda/(4M)$ and $1$ ).\\
This implies that for every $n\geq N$ and every $s\in [-\delta n^{1/2k};\delta n^{1/2k}]$,
$$n \left(G_n\left(m+sn^{-1/2k}\right)-G_n(m)\right)\geq  \frac{\lambda}{2} \frac{s^{2k}}{(2k)!}-\sum_{j=1}^{2k-1}|s|^j.$$

\begin{lem}\label{lemma4}
Let $V$ be any closed subset of $\bbR$ containing no global minima of $G$. There exists $\varepsilon>0$ such that 
$$e^{ng}\int_{V}\exp{(-nG_n(s))}\, {\rm d}s = {\cal O}\Big(e^{-n\varepsilon}\Big).$$
\end{lem}
\noindent{\bf Proof:}\\*
Since $|S_n|\leq n$, the expectation $\bbE_x(\exp(\beta J s S_n))$ is bounded above by $\exp (n\beta J |s|)$ and the function
$G_n$ satisfies for every $s\in\bbR$
$$G_n(s)\geq \frac{\beta J}{2}s^2 - \beta J |s|.$$ 
This implies that for any $s$ such that $|s|\geq 3$, $G_n(s) \geq \frac{3\beta J}{2}> g$ (since $G(0)=0$, $g$ is nonpositive).\\*
{} From Lemma \ref{lemma2}, on the compact set $W=V\cap \{s; |s|\leq 3 \}$, the sequence of functions $G_n$ converges uniformly to $G$. Hence,
the sequence $\inf_W G_n(s)$ converges to $\inf_W G(s)>g$. Let 
$$\varepsilon= \min(\frac{1}{2}(\inf_W G(s)-g),\frac{3\beta J}{2}-g) >0.$$
Then, for large $n$, $G_n(s)\geq g+\varepsilon$ on the set $V$. Hence, for large $n$,
$$e^{ng}\int_{V}\exp{(-n G_n(s)) }\, {\rm d}s \leq e^{ng} e^{-(n-1)(g+\varepsilon)} \int_{\bbR}e^{-G_n(s)}\, {\rm d}s.$$
The inequality
$$e^{-G_n(s)}\leq \exp\left(-\frac{\beta J}{2}s^2+\beta J |s|\right)$$
implies that 
$$e^{ng}\int_{V}\exp{(-n G_n(s)) }\, {\rm d}s  \leq  e^{g}e^{-(n-1)\varepsilon} \int_{\bbR}\exp\left(-\frac{\beta J}{2}s^2+\beta J |s|\right)\, {\rm d}s={\cal O}(e^{-n\varepsilon}).$$

\subsection{Proof of Theorem \ref{theo1}}

For every $y\in E$ fixed, the function $s\mapsto L(f(y),s)$ is real analytic. Moreover, if $s\in [-\alpha,\alpha]$, then for every $y\in E$, $|L(f(y),s)|\leq \alpha $. This implies that the functions 
$s\mapsto \int_E L(f(y),s) \, {\rm d}\mu (y)$ and $G$ are real analytic.\\
\\*
We now prove that $G$ has a finite number of global minima. The function $G$ goes to infinity as $s$ goes to infinity, since for every $s\in\bbR$,
 $$G(s)\geq \frac{\beta J}{2}s^2 - \beta J |s|.$$ 
This implies that the continuous function $G$ is bounded below and has a global minimum. Furthermore, the set where $G$ achieves its mimimum is bounded. Since $G$ is a non constant analytic function, the set where its first derivative $G^{(1)}$ vanishes is discrete. The set where $G$ is minimum is thus discrete. Being also bounded, it must be a finite set. 
\\*
{\bf Remark:}
It is worth remarking that every minimum point of $G$ is of finite type. Indeed, if $G$ has a minimum $m$ of infinite order, then $G^{(j)}(m)=0$ for every $j\geq 1$. Since $G$ is analytic, this implies that $G$ is constant. It contradicts the fact that $G$ tends to infinity when $s$ goes to infinity. 

\subsection{Proof of Theorem \ref{theo2}}

A consequence of Lemma \ref{lemma1} is that the density of $\frac{S_n}{n}+\frac{Y}{\sqrt n}$ is equal to
$$\frac{ \exp(-nG_n(s))}{\int \exp(-nG_n(s))\,{\rm d}s}.$$
Note that the contribution of the Gaussian random variable $\frac{Y}{\sqrt n}$ vanishes in the limit $n\rightarrow+\infty$. We thus have to prove that for any bounded continuous function $h$,
$$\frac{\int \exp(-nG_n(s))h(s)\,{\rm d}s}{\int \exp(-nG_n(s))\,{\rm d}s} = \frac{\sum_{i=1}^r b_{i,n} h(m_i)}{\sum_{i=1}^r b_{i,n}}+o(1)$$
For every $i\in\{1,\ldots,r\}$, $m_i$ is a minimum of $G$ of type $2k_i$ and strength $\lambda_i$. For every minimum $m_i$, we apply Lemma \ref{lemma3}. Then, for every $i\in\{1,\ldots,r\}$, there exist $\delta_i>0$ and $N_i\geq 1$ such that $(\ref{DL2})$ holds. Let $N$ be the maximum of $(N_i)_{i=1,\ldots,r}$ and let $\delta$ be such that $\delta\leq\min_{i=1,\ldots,r}\delta_i$ and the sets $]m_i-\delta,m_i+\delta[$ be disjoint. Let $V$ be the closed set
$$V=\bbR \setminus \bigcup_{i=1}^r \ ]m_i-\delta,m_i+\delta[.$$ 
Lemma \ref{lemma4} yields 
$$\int_V \exp(-nG_n(s))h(s)\,{\rm d}s= {\cal O}\Big(e^{-n\varepsilon-ng}\Big).$$
Use a change of variables and Lemma \ref{lemma3} to estimate the contribution of the set $]m_i-\delta,m_i+\delta[$,
\begin{eqnarray*}
 & &\int_{m_i-\delta}^{m_i+\delta}\exp(-nG_n(s))h(s)\,{\rm d}s \\ 
 &=& n^{-1/2k_i}e^{-nG_n(m_i)}\int_{-\delta n^{1/2k_i}}^{\delta n^{1/2k_i}}\exp\left(-n(G_n(m_i+sn^{-1/2k_i})-G_n(m_i))\right)h(m_i+sn^{-1/2k_i})\,{\rm d}s\\
 &=& n^{-1/2k_i}e^{-nG_n(m_i)}\left[ \int_{-\infty}^{+\infty} \exp(-\lambda_i s^{2k_i}/(2k_i)!) h(m_i) \, {\rm d}s +o(1)\right]
\end{eqnarray*}
where the last equality is obtained by combining (\ref{DL1}), (\ref{DL2}) and dominated convergence.
Hence, the integral is equivalent as $n$ goes to infinity to
$$ b_{i,n}\ h(m_i)$$
where
$$b_{i,n}= n^{-1/2k_i}e^{-nG_n(m_i)} \lambda_i^{-1/2k_i} \int_{-\infty}^{+\infty} \exp(-s^{2k_i}/(2k_i)!)\, {\rm d}s .$$
This yields the asymptotic
$$\frac{\int \exp(-nG_n(s))h(s)\,{\rm d}s}{\int \exp(-nG_n(s))\,{\rm d}s}= \displaystyle\frac{\sum_{1\leq i\leq r} b_{i,n} h(m_i)}{\sum_{1\leq i\leq r} b_{i,n}}+o(1).$$
\\*
To prove the large deviations property, we use Laplace method. From Theorem \ref{th1}, the distribution of $S_n/n$ under $\bbP_{x}$ satisfies a large deviation principle with speed $n$ and good rate function $\Lambda^\star$. The distribution $\bbQ_{n,x}$ is absolutely continuous with respect to $\bbP_x$ with density
$$\frac{{\rm d}\bbQ_{n,x}}{{\rm d}\bbP_x}=\frac{1}{\tilde{Z}_{n,x}} \exp\left( n\frac{\beta J (S_n/n)^2}{2}\right).$$
Since $S_n/n$ takes its values in $[0,1]$ and that the function $z\mapsto \beta J z^2/2$ is continuous and bounded on $[0,1]$, it follows from Varadhan integral lemma (see \cite{Dem}) that the distribution of $M_n/n$ under $\bbQ_{n,x}$ satifies a large deviation principle with speed $n$ and good rate function 
$$I_{\beta,x}(z)=  \Lambda^{*}(z) - \frac{\beta J}{2}z^2 - \inf_{z\in\bbR} \{ \Lambda^{*}(z) - \frac{\beta J}{2} z^2\}.$$

\subsection{Proof of Theorem \ref{theo3}}
{} From Lemma \ref{lemma1}, the probability density function of 
$$\frac{M_n-nm}{n^{1-1/2k}}+\frac{Y}{n^{1/2-1/2k}}$$
is given by 
$$ \frac{ \exp(-nG_n(m+sn^{-1/2k}))}{\int \exp(-nG_n(m+sn^{-1/2k}))\,{\rm d}s}.$$
The theorem is a consequence of the following result: for any bounded continuous function $h$,
\begin{equation}\label{eqres}
\frac{ \int \exp(-nG_n(m+sn^{-1/2k}))h(s)\,{\rm d}s}{\int \exp(-nG_n(m+sn^{-1/2k}))\,{\rm d}s} \rightarrow \frac { \int \exp(-\lambda s^{2k}/(2k)!)h(s)\,{\rm d}s}{\int \exp(-\lambda s^{2k}/(2k)!)\,{\rm d}s}
\end{equation}
Let $\delta$ given by Lemma \ref{lemma3} and let $V=\bbR\setminus \ ]m-\delta,m+\delta[$. Lemma \ref{lemma4} yields the following estimation 
\begin{eqnarray}\label{eqone}
 & &\int_{|s|\geq \delta n^{1/2k}} \exp(-nG_n(m+sn^{-1/2k}))h(s)\,{\rm d}s\nonumber\\
 &=& n^{1/2k}\int_V \exp(-nG_n(s))h((s-m)n^{1/2k})\,{\rm d}s\nonumber\\
 &=& {\cal O}\Big(n^{1/2k}e^{-n\varepsilon-ng}\Big)
\end{eqnarray}
From Lemma \ref{lemma3} and dominated convergence,
\begin{eqnarray}\label{eqtwo}
 & &\int_{|s|< \delta n^{1/2k}} \exp(-nG_n(m+sn^{-1/2k}))h(s)\,{\rm d}s\nonumber\\
 &=& e^{-nG_n(m)} \int_{|s|< \delta n^{1/2k}} \exp\left(-n(G_n(m+sn^{-1/2k})-G_n(m))\right)h(s)\,{\rm d}s\nonumber\\
 &=& e^{-nG_n(m)} \left[ \int_{\bbR} \exp(-\lambda s^{2k}/(2k)!)h(s)\,{\rm d}s +o(1)\right]
\end{eqnarray}
By combining equations (\ref{eqone}) and (\ref{eqtwo}) we obtain (\ref{eqres}).

\section{The phase transition}
In this section, we keep the hypotheses of the previous section and we focus on the case where $\int_E f\,{\rm d}\mu = \frac{1}{2}$ and $a=\int_E 4f(1-f)\,{\rm d}\mu>0$. We prove that a phenomenon of phase transition occurs: at high temperature (i.e. $\beta$ small), the system has the same behaviour as at infinite temperature (i.e. $\beta=0$). We prove that for $\beta<\beta_c$, the magnetization vanishes in the thermodynamical limit (i.e. $M_n/n$ converges in distribution to zero) and we also give a caracterization of the critical inverse temperature $\beta_c$.
The study of the system at lower temperature ($\beta>\beta_c$) is quite difficult since the study of the minima of $G$ strongly depends on the dynamical system $S$ and on the function $f$. We will give some general conditions under which the minima of $G$ are well known, and also study some examples. 

\subsection{The critical temperature}
In this section the parameter $\beta$ is not fixed and we write $G_{\beta}$ instead of $G$ to enhance this dependency.

\begin{defn}
Let $\beta_c$ be the critical inverse temperature defined by the relation
$$\beta_{c}=\sup\left\{\beta>0\ | \forall s\in\bbR\ ,\ G_{\beta}(s)\geq 0\right\}.$$
\end{defn}
\noindent We recall that $G_{\beta}(s)=\frac{\beta J}{2}s^2-\Lambda(\beta J s)$,
where 
$$\Lambda(y)=\int_E \log\left(f(x)e^y+(1-f(x))e^{-y}\right)\,  d\mu(x).$$
It follows from the change of variables $u=\beta J s$ that $\beta_c$ is also defined by
$$\beta_{c}=\sup\left\{\beta>0\ \Big| \forall u\in\bbR\ ,\ \Lambda(u)\leq \frac{u^2}{2\beta J}\right\}.$$
Equivalently, the critical temperature is defined by the more explicit relation
$$\frac{1}{2\beta_cJ }=\sup_{u\in\bbR^\star} \frac{\Lambda(u)}{u^2}.$$

\begin{pr}\label{criticaltemp}
The critical inverse temperature $\beta_c$ satisfies 
$$\frac{1}{J}\leq \beta_c \leq \frac{1}{Ja}.$$ 
Furthermore, if $\beta<\beta_c$, then $G_{\beta}$ has a unique minimum at point $0$.\end{pr}
\noindent{\bf Proof :}\\
We use the fact that $\Lambda(u)\sim \frac{a}{2}u^2$ as $u\rightarrow 0$. Hence, $u\mapsto\frac{\Lambda(u)}{u^2}$ is continuous on $\bbR^\star$, tends to $a/2$ as $u\rightarrow 0$ and to $0$ as $u\rightarrow \pm\infty.$ This implies that 
$$\frac{a}{2} \leq \sup_{u\in\bbR^\star} \frac{\Lambda(u)}{u^2} $$
and then $\beta_c \leq \frac{1}{Ja}$. 
Note that the equality $\beta_c = \frac{1}{Ja}$ holds if and only if $u\mapsto \frac{\Lambda(u)}{u^2}$ reaches its maximum at point $0$. The inequality $\beta_c \geq \frac{1}{J}$ follows from the fact that the second derivative of $\Lambda$ is bounded above by $1$ which implies that $\Lambda(u)\leq \frac{1}{2}u^2. $ \\
The relation 
$$G_{\beta}\left(\frac{s}{\beta J}\right)=G_{\beta_c}\left(\frac{s}{\beta_c J}\right)+\left(\frac{1}{2\beta J}-\frac{1}{2\beta_c J}\right) s^2$$
implies that if $\beta < \beta_c$, the function $G_{\beta}$ has a unique minimum at point $0$ equal to $0$.

\begin{theo}\label{theo4}
Let us assume that $f\in\cC_{1/2}(S)$. Then, for every $\beta<\beta_c$, the distribution of $M_n$ under $\bbQ_{n,x}$ satisfies a law of large numbers:
$$\frac{1}{n} M_n \Rightarrow \delta_0\ ,\ {\rm\ as\ } n\rightarrow\infty $$
and a central limit theorem:  
$$\frac{1}{\sqrt n} M_n \Rightarrow \cN(0,\sigma^2)\ ,\ {\rm\ as\ } n\rightarrow\infty$$
with $\sigma^2=a/(1-\beta J a).$
\end{theo}
\noindent{\bf Proof:}\\
The law of large numbers is an application of Theorem \ref{theo2} and the central limit theorem is an application of Theorem \ref{theo3}. We verify that the assumptions of these theorems are satisfied. From Proposition \ref{criticaltemp}, the hypothesis $\int_{E}f\, d\mu=\frac{1}{2}$  and $\beta<\beta_c$ implies that the function $G_{\beta}$ has a unique minimum at point $m=0$ of type $2k=2$ and strength  $\lambda=G_{\beta}^{(2)}(0)=\beta J ( 1-\beta J a)>0$.
Since the function $y\mapsto \frac{\partial }{\partial s}L(f(y), 0)$ is equal to  $\beta J(2f-1)$, the hypothesis $f\in\cC_{1/2}(S)$ implies that $y\mapsto \frac{\partial }{\partial s}L(f(y), 0)$ belongs to $\cC_{1/2}(S)$. Moreover, the function $y\mapsto \frac{\partial^2 }{\partial s^2}L(f(y), 0)$ is equal to $(\beta J)^2\left(4f(1-f)\right)$ which belongs to $\cC_{1}(S)$ since $f$ is continuous.

\subsection{General study for a specific class of systems}
We consider here a class of systems for which we can study the minima of the function $G_{\beta}$ for every $\beta>0$. 
We suppose that the system satisfies the following hypothesis 
$${\bf (H)} \hspace{0.5cm} {\rm The\  function\ } \Lambda {\rm \  is\ even\ and\ its \ derivative \ is\ concave\ on\ } (0,+\infty).    $$
Before introducing our main results, we exhibit some cases where the hypothesis $(H)$ is satisfied. Note that $(H)$ is an assumption on the function 
$$\Lambda:u\mapsto \int_E \log\left(f(x)e^u+(1-f(x))e^{-u}\right)\,{\rm d}\mu(x)$$
that does not really depend on the dynamical system $S=(E,\cA,\mu,T)$ but only on the image distribution of $\mu$ under the application $f$, that we denote by $\mu_f $ that is the measure on $[0,1]$ such that for every Borel set $A$, $$\mu_f(A)=\mu(f^{-1}(A)).$$  

\begin{pr}
The hypothesis $(H)$ is satisfied in the following cases: \vspace{-0.2cm}
\begin{enumerate}
        \item the measure $\mu_f $ is equal to ${\bf 1}_{[0,1]}(x)\,{\rm d}x$.
        \item the measure $\mu_f $ is equal to $\frac{1}{2}(\delta_{\lambda}+\delta_{1-\lambda})$, with 
        $ \frac{1}{2}-\frac{\sqrt 3}{6}\leq \lambda\leq \frac{1}{2}$, (and if $0\leq\lambda<\frac{1}{2}-\frac{\sqrt 3}{6}$, the hypothesis $(H)$ is not satisfied.)
        \item the measure $\mu_f $ has its support included into $[\frac{1}{2}-\frac{\sqrt 3}{6},\frac{1}{2}+\frac{\sqrt 3}{6}]$  and satisfies the symmetry condition: $\mu_f =\mu_{1-f} $.
\end{enumerate}
\end{pr}
\noindent{\bf Proof:} \\
{\it 1.:\ }In this case, the function $\Lambda$ is equal to
$$\Lambda(u)=\int_0^1 \log\left(xe^u+(1-x)e^{-u}\right)\,{\rm d}x
 =\left\{\begin{array}{ll}
 \frac{u}{\tanh(u)}-1 & {\mbox if}\ \    u\neq 0 \\
 \  \  \ 0  &  {\mbox if}\ \   u=0. 
 \end{array}\right. 
  $$
It is an even function. In order to prove the concavity of $\Lambda'$, we compute the third derivative $\Lambda^{(3)}$ for $u\in (0,+\infty)$ :
$$\Lambda^{(3)}(u)=  \frac{2(1-\tanh(u)^2)}{\tanh(u)^4}(3\tanh(u)-3u+u\tanh(u)^2)$$
which is negative on $(0,+\infty)$, so $\Lambda$ satisfies the hypothesis $(H)$.\\
{\it 2.:\ } When $\mu_f=\frac{1}{2}(\delta_{\lambda}+\delta_{1-\lambda})$, the function $\Lambda$ is equal to
\begin{eqnarray*}
\Lambda(u)&=&\frac{1}{2}\left[\log(\lambda\ e^u+(1-\lambda)\ e^{-u})+\log((1-\lambda)\ e^u+\lambda\ e^{-u})\right]\\
 &=& \log(\cosh(u))+\frac{1}{2}\log\left(1-(2\lambda-1)^2\tanh(u)^2\right),
\end{eqnarray*}
which is an even function and its third derivative is given for $u\in (0,+\infty)$ by
$$\Lambda^{(3)}(u)=\frac{2t(1-t^2)(1-g^2)}{\left(1-g^2t^2\right)^3}\left[g^2(g^2-3)t^2+(3g^2-1)\right]$$
with $t=\tanh(u)\in(0,1)$ and $g=2\lambda-1\in\ (-1,1)$.
The fraction is nonnegative and for every $\lambda\in [\frac{1}{2}-\frac{\sqrt 3}{6},\frac{1}{2}]$, the bracket is nonpositive since  $g^2(g^2-3)t^2+(3g^2-1)\leq (3g^2-1)\leq 0$. Thus, $\Lambda^{(3)}$ is negative on $(0,+\infty)$ and assumption $(H)$ is satisfied. \\
{\it 3.:\ } The condition that  $f$ and $1-f$ have the same distribution under $\mu$ implies that the function $\Lambda$ is even : for every $u\in\bbR$,
\begin{eqnarray*}
\Lambda(u)&=&\int_E \log\left(f(x)e^u+(1-f(x))e^{-u}\right)\,{\rm d}\mu(x)\\
 &=&\int_E \log\left((1-f(x))e^u+(1-(1-f(x)))e^{-u}\right)\,{\rm d}\mu(x)\\
 &=&\Lambda(-u).
\end{eqnarray*}
Writing $\Lambda(u)=\frac{\Lambda(u)+\Lambda(-u)}{2}$ yields the following expression :
$$\Lambda(u)=\log(\cosh(u))+\frac{1}{2}\int_E \log\left(1-(2f(x)-1)^2\tanh(u)^2\right)\,{\rm d}\mu(x).$$
We compute the third derivative using this last formula. This yields
$$\Lambda^{(3)}(u)=2t(1-t^2)\int_E \frac{1-g(x)^2}{\left(1-g(x)^2t^2\right)^3}\left[g(x)^2(g(x)^2-3)t^2+(3g(x)^2-1)\right]\,{\rm d}\mu(x)$$
with $t=\tanh(u)\in(-1,1)$ and $g=2f-1$. When  the measure $\mu_f$ has its support included into $[\frac{1}{2}-\frac{\sqrt 3}{6},\frac{1}{2}+\frac{\sqrt 3}{6}]$, then $(3g^2-1)$ is $\mu$-almost everywhere nonpositive and hence  $\Lambda^{(3)}(u)$ and $-u$ are of the same sign, so assumption $(H)$ is satisfied.

We are now ready to give a complete description of the minima of the function $G_{\beta}$ in function of the inverse temperature $\beta$ and to show a phase transition at the critical inverse temperature $\beta=\beta_c$.
\begin{theo}
Let us assume that hypothesis $(H)$ is satisfied. 
Then, the following results hold:
\begin{enumerate}
\item the critical inverse temperature $\beta_c$ is equal to $1/(Ja)$.
\item for $\beta<\beta_c$, $G_{\beta}$ admits $0$ as unique minimum of type $2$ and strength $\beta J(1-\beta J a)$.
\item for $\beta=\beta_c$, $G_{\beta}$ admits $0$ as unique minimum of type $\geq  4$. 
\item for $\beta>\beta_c$, $G_{\beta}$ admits two global minima at points $\pm m$ (with $m>0)$, of same type equal to 2 and same strength. The point $m$ called {\it spontaneous magnetization} is the unique positive solution
of the equation
$$ m=\int_E \frac{\tanh(\beta Jm) +(2f-1)}{1+(2f-1)\tanh(\beta Jm)}\ {\rm d} \mu.$$
\end{enumerate}
\end{theo}
\noindent{\bf Proof:}\\
{\it 1.:\ } We prove that for every $u\in\bbR$, 
\begin{equation}\label{ineq}
\Lambda(u)\leq \frac{a}{2}u^2.
\end{equation}
This implies that 
$$\frac{1}{2\beta_c J}=\sup_{u\in\bbR^\star} \frac{\Lambda(u)}{u^2}\leq \frac{a}{2}$$
and hence that $\beta_c\geq \frac{1}{Ja}$. From Proposition $\ref{criticaltemp}$, the equality is proved.\\
We now prove inequality $(\ref{ineq})$ using hypothesis $(H)$.  Since $\Lambda'$ is concave on $(0,+\infty)$, its derivative $\Lambda^{(2)}$ is a nonincreasing function. Hence, for every $u>0$, $\Lambda^{(2)}(u)\leq \Lambda^{(2)}(0)=a $. The function $\Lambda'$ vanishes at point $0$ and has a derivative bounded above by $a$: this implies that $\Lambda'(u)\leq au$. Integrating one more time yields $\Lambda(u)\leq \frac{a}{2}u^2$ for every $u>0$. Since the function $\Lambda$ is even, the same inequality holds for $u<0$ and this proves inequality $(\ref{ineq})$.\\
{\it 2.:\ } This result is contained in Proposition \ref{criticaltemp}.\\
{\it 3.:\ } From assumption $(H)$ the function $G_{\beta}'(s)$ equal to $\beta J (s -\Lambda'(\beta J s))$ is odd and convex on $(0,+\infty)$. For $\beta=\beta_c$, $G_{\beta}^{'}(0)=G_{\beta}^{(2)}(0)=0$, so from convexity, $G_{\beta}^{'}$ is nonnegative on $(0,+\infty)$. The point $0$ is the only point where $G_{\beta}^{'}$ vanishes. Otherwise, if there exists some $u>0$ such that $G_{\beta}^{'}(u)=0$, the convex function would be identically zero on $[0,u]$ , and being real analytic, it would be identically zero on $\bbR$, which is not the case. Thus zero is the only point where $G_{\beta}^{'}$ vanishes, and the function $G_{\beta}$ has a unique global minimum at point $0$.\\
{\it 4.:\ } For $\beta>\beta_c$, under assumption $(H)$, the function $G_{\beta}^{'}$ is odd and convex on $(0,+\infty)$, moreover, $G_{\beta}^{'}(0)=0$, $G_{\beta}^{(2)}(0)=\beta J(1-\beta J a)<0$ and $\lim_{s\rightarrow +\infty}G_{\beta}^{'}(s)=+\infty$. Hence, there exists an unique real $m\in (0,+\infty)$ such that $G_{\beta}^{'}(m)=0$, thus the function $G_{\beta}^{'}$ vanishes only at points $-m$, $0$ and $m$. The function $G_{\beta}$ reaches its global minimum at points $-m$ and $m$, and has a local maximum at $0$. Since $\Lambda$ is even, the minima $m$ and $-m$ have same type and same strength. \\*
We now prove that the type of $m$ is equal to 2: by the mean value theorem there exists $m_0\in (0, \beta Jm)$ such that
$$ \Lambda^{(2)}(m_0) = \frac{\Lambda^{'}(\beta Jm)}{\beta Jm}=  \frac{1}{\beta J}.$$
The real analytic function $\Lambda^{(2)}$ is not constant on $[0,+\infty[$ (since $\Lambda^{(2)}(0)=a>0$ and $\lim_{s\rightarrow +\infty} \Lambda^{(2)}(s)=0$). Then,  since $\Lambda^{(3)}$ is nonpositive on $(0,+\infty)$, the function $\Lambda^{(2)}$ is strictly decreasing on $(0,+\infty)$. It follows that $\Lambda^{(2)}(m_0)> \Lambda^{(2)}(\beta Jm)$ and that
$$  G_{\beta}^{(2)}(\beta Jm)= (\beta J)^{2} \left( \frac{1}{\beta J} -\Lambda^{(2)}(\beta Jm)\right)= (\beta J)^{2} \left( \lambda^{(2)}(m_0) -\Lambda^{(2)}(\beta Jm)\right) >0.$$
This completes the proof of the theorem.

As a consequence, the asymptotic behaviour of $M_n$ under $\bbQ_{n,x}$ for $\beta\geq \beta_c$ is deduced for the systems satisfying assumption $(H)$. Recall that Theorem $\ref{theo4}$ treats the case $\beta<\beta_c$ for general systems. 
\begin{theo}\label{theo5}
Assume that assumption $(H)$ is satisfied. 
\begin{enumerate}
        \item  When $\beta=\beta_c$, denote by $2k\geq 4$ and $\lambda$ the type and the strength of $m=0$ the minimum of $G_{\beta}$. Assume that for every $j\in\{1,\ldots,2k\}$, the function $f^j$ belongs to $\mathcal{C}_{j/2k}(S).$
Then, $$\frac{M_n}{n} \Rightarrow \delta_0$$
and
$$\frac{M_n}{n^{1-1/2k}}\Rightarrow Z $$
where $Z$ is the probability measure with density function
$$C \exp\left(-\lambda s^{2k}/(2k)!\right),$$
$C$ being the normalizing constant.
        \item When $\beta> \beta_c$, assume that the functions $\frac{\partial}{\partial s}L(f(.),\pm \beta J m)$ belong to the set $\mathcal{C}_{1/2}(S).$ \\*
Then, for every bounded continuous function $h$, the expectation  of $h(M_n/n)$ under $\bbQ_{n,x}$ is equivalent, as $n$ goes to infinity, to
$$\frac{b_{m,n}\ h(m) +b_{-m,n}\  h(-m)}{b_{m,n}+b_{-m,n}}$$
where 
$$b_{m,n}=\prod_{j=1}^n\left(f(T^jx)\ e^{\beta Jm} +(1-f(T^jx))\ e^{-\beta Jm}\right).$$ 
   
\end{enumerate}
\end{theo}
\noindent{\bf Remark:} A straightforward computation gives 
$$G_{\beta}^{(4)}(0)=2(\beta J)^4 (3I_4-4I_2+1)$$ 
where $I_2=\int_E (2f-1)^2\ {\rm d}\mu$ and $I_4=\int_E (2f-1)^4\ {\rm d}\mu$.
So, when $\beta=\beta_c$, the type of 0 is equal to 4 if and only if $3I_4-4I_2+1 >0$. It is always verified when the support of $f$ is strictly included into $[\frac{1}{2}-\frac{\sqrt 3}{6},\frac{1}{2}+\frac{\sqrt 3}{6}]$.   

\section{A particular case: the rotation on the torus} 
In this section we treat the particular case of the irrational rotation on the torus which corresponds to a quasiperiodic random field already mentioned in \cite{KPZ}. 
It is one of the dynamical systems $S$ in ergodic theory for which we are able to provide a large subclass of ${\cal C}_{\alpha}(S)$. In the first section we 
first give a precise description of this subclass in terms of the diophantine properties of the irrational angle. In the second one we apply results of Section 4 to this particular dynamical system when the function $f$ is the identity function.    
\subsection{Some Results on Diophantine Approximations}
Let us consider the dynamical system $(\bbT^{r},\cB(\bbT^{r}),\lambda,T_{\alpha})$ where $\lambda$ is the Lebesgue measure on the torus $\bbT^{r}$ and $T_{\alpha}$ is the irrational rotation over $\bbT^{r}$ defined by $x\ra x+\alpha \mod 1$. It is well known that under these conditions this dynamical system is ergodic and for every $f\in L^{1}(\lambda)$, for almost every $x\in\bbT^{r}$,
$$M_{n}=\frac{1}{n}\sum_{l=1}^{n}f(T_{\alpha}^{l}x)-\int_{\bbT^{r}}f(t)dt\ra_{n\ra\infty} 0$$
When $f$ is with bounded variation, this result holds for every $x\in\bbT^{r}$ and it is possible to determine the speed of convergence of the sequence $M_{n}$ to 0 in terms of arithmetic properties of the irrational vector $\alpha$. When $r=1$, for all irrational badly approximated by rationals, Denjoy-Koksma's inequality gives us a majorization of $M_{n}$ uniformly in $x$ for $n$ large enough. But when $r\geq 2$, Denjoy-Koksma's inequality does not hold (see Yoccoz \cite{yo}) and the method of low discrepancy sequences has to be used.
\\
\\
\subsubsection{Case of one-dimensional torus}
Let $\alpha$ be an irrational.\index{irrational angle}
We call a rational $\frac{p}{q}$ with $p, q$ relatively prime such that $|\alpha-\frac{p}{q}|<
\frac{1}{q^2}$, a rational approximation of $\alpha$.
When $\alpha$ has the continued fraction expansion $\alpha=[\alpha]+[a_{1},
\ldots,a_{n},\ldots]$, the $n$-th principal convergent of $\alpha$ is $\frac{p_{n}}{q_{n}}$
where, $\forall n\geq 2$,
$$p_{n}=a_{n}p_{n-1}+p_{n-2}$$
$$q_{n}=a_{n}q_{n-1}+q_{n-2};$$
the recurrence is given by defining the values of $p_{0}, p_{1}$ and $q_{0},q_{1}$.
\\
\\
{\bf Denjoy-Koksma's inequality}\index{Denjoy-Koksma's inequality}
{\it Let $f:\bbR\rightarrow [0,1]$ be a function with bounded variation $V(f)$ and $\frac{p}{q}$ a
rational approximation of $\alpha$.
Then, for every $x\in \bbT^{1}$,
$$|\sum_{l=1}^{q}f(T_{\alpha}^{l}x)-q\int_{\bbT^1}f(t)dt|\leq V(f).$$
}\\

\begin{pr}
Let $f$ be a function with bounded variation $V(f)$. For every
 irrational $\alpha$ such that the inequality $a_{m}<m^{1+\epsilon},$
where $\epsilon>0,$ is satisfied eventually for all $m$,
$$\sup_{x\in\bbT^1}\Big|\sum_{l=1}^{n}\big(f(T_{\alpha}^{l}x)-\int_{\bbT^{1}}f(t)dt\big)\Big|=\cO(\log^{2+\epsilon} n).$$
\end{pr}
\noindent{\bf Proof:}\\*
The sequence of integers $(q_{i})_{i\geq 1}$ being strictly increasing,
for a given $n\geq 1$, there exists $m_{n}\geq 0$ such that
$$q_{m_{n}}\leq n<q_{m_{n}+1}.$$
By Euclidean division, we have $n=b_{m_{n}}q_{m_{n}}+n_{m_{n}-1}$
with $0\leq n_{m_{n}-1}<q_{m_{n}}$.
We can use the usual relations
$$q_{0}=1, q_{1}=a_{1}$$
\begin{equation}\label{eq}
q_{n}=a_{n}q_{n-1}+q_{n-2}, n\geq 2.
\end{equation}
We obtain that $(a_{m_{n}+1}+1)q_{m_{n}}>q_{m_{n}+1}>n$ and so
 $b_{m_{n}}\leq a_{m_{n}+1}$.
If $m_{n}>0$, we may write $n_{m_{n}-1}=b_{m_{n}-1}q_{m_{n}-1}+n_{m_{n}-2}$
with $0\leq n_{m_{n}-2}<q_{m_{n}-1}.$
Again, we find $b_{m_{n}-1}\leq a_{m_{n}}$.
Continuing in this manner, we arrive at a representation for $n$ of the form
$$n=\sum_{i=0}^{m_{n}}b_{i}q_{i}$$
with $0\leq b_{i}\leq a_{i+1}$ for $0\leq i\leq m_{n}$ and $b_{m_{n}}\geq 1$.
Using Denjoy-Koksma's inequality, we get
\begin{eqnarray*}
|\sum_{l=1}^{n}f(T_{\alpha}^{l}x)-n\int_{\bbT^1}f(x)dx|&\leq&
V(f)\sum_{i=0}^{m_{n}}b_{i}\\
&\leq& V(f)\sum_{i=0}^{m_{n}}a_{i+1}.
\end{eqnarray*}
By hypothesis, there exists $m_{0}\geq 1$ such that,
$$a_{m}<m^{1+\epsilon}, \forall m\geq m_{0}.$$
Let $n$ be such that $m_{n}>m_{0}$.
Thus,
$$|\sum_{l=1}^{n}f(T_{\alpha}^{l}x)-n\int_{\bbT^1}f(t)dt|\leq
V(f)(\sum_{i=0}^{m_{0}-1}a_{i+1}+(m_{n}+1)^{2+\epsilon}).$$ We
need to know the asymptotic behavior of $m_{n}$. When $\alpha$ is
the golden ratio, $a_{n}=1,\ \forall n\geq 1$ and the relation
(\ref{eq})
 implies that $q_{n}\sim \frac{1}{\sqrt{5}}\alpha^{n+1}$.
Let $\alpha'$ be another irrational; its partial quotients $a_{n}'$ satisfy necessarily
$a_{n}'\geq1$.
Using the relation (\ref{eq}), we see that $q_{n}'\geq q_{n}, \forall n\geq 1$.
Therefore, $m_{n}=\cO(\log n)$ and the proposition is proved.

\subsubsection{Generalization to $r-$dimensional torus}
We recall some definitions and well known results from the method of low
discrepancy sequences in dimension $r\geq 1$.

Suppose we are given a function $f(x)=f(x^{(1)},\ldots,x^{(r)})$ with $r\geq 1.$
By a partition $P$ of $[0,1]^r$, we mean a set of $r$ finite sequences $\eta_{0}^{(j)},
\eta_{1}^{(j)},\ldots,\\
\eta_{m_{j}}^{(j)}(j=1,\ldots,r),$ with $0=\eta_{0}^{(j)}\leq
\eta_{1}^{(j)}\leq\ldots\leq\eta_{m_{j}}^{(j)}=1$
for $j=1,\ldots,r$.
In connection with such a partition, we define, for $j=1,\ldots,r$ an operator
$\Delta_{j}$ by
$$\Delta_{j}f(x^{(1)},\ldots,x^{(j-1)},\eta_{i}^{(j)},x^{(j+1)},\ldots,x^{(r)})
=f(x^{(1)},\ldots,x^{(j-1)},\eta_{i+1}^{(j)},$$
$$x^{(j+1)},\ldots,x^{(r)})
-f(x^{(1)},\ldots,x^{(j-1)},\eta_{i}^{(j)},x^{(j+1)},\ldots,x^{(r)}),$$
for $0\leq i< m_{j}.$
\begin{defn}
\begin{enumerate}
\item For a function $f$ on $[0,1]^r$, we set
$$V^{(r)}(f)=\sup_{P}\sum_{i_{1}=0}^{m_{1}-1}\ldots \sum_{i_{r}=0}^{m_{r}-1}
|\Delta_{1,\ldots,r}f(\eta_{i_{1}}^{(1)},\ldots,\eta_{i_{r}}^{(r)})|,$$
where the supremum is extended over all partitions $P$ of $[0,1]^r$.
If $V^{(r)}(f)$ is finite, then $f$ is said to be {\it of bounded variation on $[0,1]^r$
in the sense of Vitali}.
\item For $1\leq p\leq r$ and $1\leq i_{1}<i_{2}<\ldots<i_{p}\leq r$,
we denote by $V^{(p)}(f;i_{1},\ldots,i_{p})$ the $p$-dimensional variation in the
sense of Vitali of the restriction of $f$ to
\begin{center}
$E_{i_{1}\ldots i_{p}}^r=\{(t_{1},\ldots,t_{r})
\in [0,1]^r; t_{j}=1$ whenever $j$ is none of the $i_{r}, 1\leq r\leq p\}.$
\end{center}
If all the variations $V^{(p)}(f;i_{1},\ldots,i_{p})$ are finite, the function $f$ is said
to be {\it of bounded variation on $[0,1]^r$ in the sense of Hardy and Krause}.
\end{enumerate}
\end{defn}
Let $x_{1},\ldots,x_{n}$ be a finite sequence of points in $[0,1]^r$ with
$x_{l}=(x_{l_{1}},\ldots,x_{l_{r}})$ for $1\leq l\leq n$.
We introduce the function
$$R_{n}(t_{1},\ldots,t_{r})=\frac{A(t_{1},\ldots,t_{r};n)}{n}-t_{1}\ldots t_{r}$$
for $(t_{1},\ldots,t_{r})\in [0,1]^r$, where $A(t_{1},\ldots,t_{r};n)$ denotes the number
of elements $x_{l}, 1\leq l\leq n,$ for which $x_{l_{i}}<t_{i}$ for $1\leq i\leq r$.
\begin{defn}
The {\it discrepancy} $D_{n}^{*}$ of the sequence $x_{1},\ldots,x_{n}$ in $[0,1]^r$ is defined
to be
$$D_{n}^{*}=\sup_{(t_{1},\ldots,t_{r})\in [0,1]^r}|R_{n}(t_{1},\ldots,t_{r})|.$$
\end{defn}
For a real number $t,$ let $\|t\|$ denote its distance to the nearest integer, namely,
\begin{eqnarray*}
\| t\|&=&\inf_{n\in\bbZ}\mid t-n\mid\\
&=&\inf(\{t\},1-\{t\})
\end{eqnarray*}
 where $\{t\}$ is the fractional part of $t$.
\begin{defn}
For a real number $\eta$, a $r$-tuple $\alpha=(\alpha_{1},\ldots,\alpha_{r})$
of irrationals is said to be of {\it type} $\eta$ if $\eta$ is the infimum of all numbers
$\sigma$ for which there exists a positive constant $c=c(\sigma; \alpha_{1},
\ldots,\alpha_{r})$ such that
$$r^{\sigma}(h)\|<h,\alpha>\|\geq c$$
holds for all $h\neq 0$ in $\bbZ^{r}$, where $r(h)=\prod_{i=1}^{r}\max(1,|h_{i}|)$ and
$<\cdot,\cdot>$ denotes the standard inner product in $\bbR^{r}.$
\end{defn}
%\N
The type $\eta$ of $\alpha$ is also equal to
$$\sup\{\gamma:\inf_{h\in(\bbZ^{r})^{*}}r^{\gamma}(h)\|<h,\alpha>\|=0\}.$$

We always have $\eta\geq 1$ (see \cite{nie}). Now we give a
result (see \cite{kui}) which yields the asymptotic behavior of
the discrepancy of the sequence
$w=(x_{1}+l\alpha_{1},\ldots,x_{r}+l\alpha_{r}), l=1,2,\ldots $ as
a function of the mutual irrationality of the components of
$\alpha$.
\begin{pr}
Let $\alpha=(\alpha_{1},\ldots,\alpha_{r})$ be an irrational vector.
Suppose there exists $\eta\geq 1$ and $c>0$ such that
$$r^{\eta}(h)\|<h,\alpha>\|\geq c$$
for all $h\neq 0$ in $\bbZ^{r}$.
Then, for every $x\in [0,1]^r$, the discrepancy of the sequence $w=(x_{1}+l\alpha_{1},
\ldots,x_{r}+l\alpha_{r}), l=1,2,\ldots $ satisfies $D_{n}^{*}(w)=\cO(n^{-1}\log^{r+1} n)$
for $\eta=1$ and $D_{n}^{*}(w)=\cO(n^{-\frac{1}{((\eta-1)r+1)}}\log n)$ for $\eta>1$.
\end{pr}
The proof is based on the Erd\"{o}s-Tur\'{a}n-Koksma's theorem:
For $h\in\bbZ^r,$ define $p(h)=\max_{1\leq j\leq r}|h_{j}|$.
Let $x_{1},\ldots,x_{n}$ be a finite sequence of points in $\bbR^{r}$.
Then, for any positive integer $m$, we have
$$D_{n}^{*}\leq C_{r}\left(\frac{1}{m}+\sum_{0\leq p(h)\leq m}\frac{1}{r(h)}\left|\frac{1}{n}
\sum_{l=1}^{n}e^{2\pi i<h,x_{l}>}\right|\right)$$
where $C_{r}$ only depends on the dimension $r$.
This theorem combined with the results of \cite{kui} (p.131)
 gives us the result.
\begin{theo}[Hlawka, Zaremba]
Let $f$ be of bounded variation on $[0,1]^r$ in the sense of Hardy and Krause, and let
$\omega$ be a finite sequence of points $x_{1},\ldots,x_{n}$ in $[0,1]^r$.
Then, we have
$$|\frac{1}{n}\sum_{l=1}^{n}f(x_{l})-\int_{\bbT^r}f(t)dt|\leq\sum_{p=1}^{r}\sum_{
1\leq i_{1}<i_{2}<\ldots<i_{p}\leq r}V^{(p)}(f;i_{1},\ldots,i_{p})
D_{n}^{*}(\omega_{i_{1}\ldots i_{p}}),$$
where $D_{n}^{*}(\omega_{i_{1}\ldots i_{p}})$ is the discrepancy in $E_{i_{1}\ldots i_{p}}^r$
of the sequence $\omega_{i_{1}\ldots i_{p}}$ obtained by projecting $\omega$
onto $E_{i_{1}\ldots i_{p}}^r$.
\end{theo}
\begin{pr}\label{pr18}
Let $f$ be a function with bounded variation in the sense of Hardy
and Krause, and $\alpha$ an irrational vector of type $\eta$, then
$$\sup_{x\in\bbT^{r}}\Big|\sum_{l=1}^{n}\big(f(T_{\alpha}^lx)-\int_{\bbT^{r}}f(t)dt\big)\Big|=
\left\{\begin{array}{ll}
\cO(\log^{r+1} n)\  \mbox{if}\ \eta=1\\
\cO(n^{1-\frac{1}{((\eta-1)r+1)}}\log n)\
\mbox{if}\  \eta>1.\\
\end{array}
\right.$$
\end{pr}
\noindent{\bf Proof:}\\*
Let $\eta'$ be such that $\eta\leq\eta'<1+\frac{1}{r}$.
There exists $c>0$ such that
$$r^{\eta'}(h)\|<h,\alpha>\|\geq c$$
holds for all $h\neq 0$ in $\bbZ^{r}$.
Suppose we are given a $p$-tuple $\alpha_{p}=(\alpha_{i_{1}},\ldots,\alpha_{i_{p}}), 1\leq p
\leq r,$ of $\alpha$, then
$$r^{\eta'}(h)\|<h,\alpha_{p}>\|\geq c$$
holds for all $h\neq 0$ in $\bbZ^{p}, 1\leq p\leq r$.
Thus, every $p$-tuple, $1\leq p\leq r$, is of type $\delta$ such that
 $1\leq \delta \leq \eta$ and $(\alpha_{i_{1}},\ldots,\alpha_{i_{p}})$
is an irrational vector.
For every $p, 1\leq p\leq r$, we define $w_{i_{1}\ldots i_{p}}$ by the projection
of $w$ on $E_{i_{1}\ldots i_{p}}^r$. From the previous proposition, we have for every $p, 1\leq p\leq r$,
$$\left\{\begin{array}{lcl}
 nD_{n}^{*}(w_{i_{1}\ldots i_{p}})&=&\cO(\log^{p+1} n)\  \mbox{if}\ \delta=1\\
 nD_{n}^{*}(w_{i_{1}\ldots i_{p}})&=&\cO(n^{1-\frac{1}{((\delta-1)p+1)}}\log n)\  \mbox{if}\ 1
<\delta \leq \eta.
\end{array}
\right.$$
Now, $\forall p=1,\ldots,r,$
$$0\leq1-\frac{1}{(\delta-1)p+1}\leq 1-\frac{1}{(\eta-1)r+1}\leq 1.$$
Therefore, using Hlawka-Zaremba's theorem, we obtain Proposition \ref{pr18}.
\subsection{A particular example: $f(x)=x$}
Consider the irrational rotation on the one-dimensional torus with angle of type $\eta$ and 
we choose $f(x)=x$. Clearly, the integral of $f$ is equal to $1/2$ and $a=2/3$.
We apply Theorems  4.2 and 4.3 for this particular example.  
\begin{theo}\label{exam}
The following results hold:
\begin{enumerate}
\item The critical inverse temperature $\beta_c$ is equal to $3/(2J)$.
\item When $\beta<\beta_c$, if $\eta<2$,
$$\frac{M_n}{n} \Rightarrow \delta_0\ ,\ {\rm\ as\ } n\rightarrow\infty $$
and   
$$\frac{M_n}{\sqrt n} \Rightarrow \cN(0,\sigma^2)\ ,\ {\rm\ as\ } n\rightarrow\infty$$
with $\sigma^2=2/(3-2 \beta J)$  .
\item When $\beta=\beta_c$, if $\eta<4/3$,
$$\frac{M_n}{n} \Rightarrow \delta_0\ ,\ {\rm\ as\ } n\rightarrow\infty$$
and
$$\frac{M_n}{n^{3/4}}\Rightarrow Z \ ,\ {\rm\ as\ } n\rightarrow\infty$$
where $Z$ is the probability measure with density function
$$\displaystyle\frac{\sqrt{3}\ \Gamma(\frac{3}{4})} {\sqrt{2}\ \pi\ \sqrt[4]{5}}   \exp\left(-9 s^{4}/80\right).$$
\end{enumerate}
\end{theo}
\noindent{\bf Proof:}\\*
Assertion 1. comes from a direct application of Theorem 4.1. To prove 2. remark that this particular example corresponds to 1. from Proposition 5.3, so $(H)$ is satisfied. Moreover,  from the remark following Theorem 4.3, it is easy to prove that the type of the unique minimum 0 is equal to 4. Finally, by combining Theorem 4.3 and Proposition \ref{pr18} we get the result.


\begin{thebibliography}{99}
\bibitem{bil} {\sc Billingsley, P.} {\it Convergence of probability measures.} {\rm Wiley, New York
    (1968).}

\bibitem{pel} {\sc Boldrighini, C.,} {\sc Minlos, R. A.,} and {\sc Pellegrinotti, A.} {\rm Almost-sure central limit theorem for a Markov model of random walk in dynamical random environment.} {\it Probab. Theory Related Fields} (1997), Vol. {\bf 109}, No 2, 245--273.

\bibitem{com} {\sc Comets, F.,} {\sc Gantert, N.,} and {\sc Zeitouni, O.} {\rm Quenched, annealed and functional large deviations for one-dimensional random walk in random environment.} {\it Probab. Theory Related Fields} (2000), Vol. {\bf 118}, No 1, 65--114.

\bibitem{Dem} {\sc Dembo, A.} and {\sc Zeitouni, O.} {\it Large Deviations Techniques and Applications.} {\rm Springer, (1998)}.

\bibitem{din} {\sc Dinwoodie, I. H.} and {\sc Zabell, S. L.} {\rm Large deviations for exchangeable random vectors.}
{\it Ann. Probab.} (1992), Vol. {\bf 20}, No 3, 1147--1166.

\bibitem{DGPSP} {\sc Dombry, C., Guillotin-Plantard, N., Pincon, B.} and {\sc Schott, R.}  {\rm Data Structures with Dynamical Random Transitions.} {\it Random Structures and Algorithms} (2006), Vol. {\bf 28}, No 4, 403 -- 426.

\bibitem{livEll} {\sc Ellis, R.S.} {\it Entropy, Large Deviations and Statistical Mechanics.} {\rm New-York, Springer-Verlag, (1985)} 

\bibitem{ell1} {\sc Ellis, R.S., Newman, C. M.}  {\rm Limit theorems for sums of dependent random variables occuring in statistical mechanics.} {\it Z. Wahrscheinlichkeitstheorie verw. Gebiete} (1978), Vol. {\bf 44}, 117--139.

\bibitem{ell2} {\sc Ellis, R.S., Newman, C. M.} and {\sc Rosen, J. S.} {\rm Limit theorems for sums of dependent random variables occuring in statistical mechanics.} {\it Z. Wahrscheinlichkeitstheorie verw. Gebiete} (1980), Vol. {\bf 51}, 153--169.

\bibitem{Fel}{\sc Feller, W.}
  {\it Introduction to Probability Theory and its Applications.} {Vol. II, \rm Wiley, New York (1971).}
  
\bibitem{Picco} {\sc Fontes, L.R.},  {\sc Mathieu, P.} and {\sc Picco, P.} {\rm On the averaged dynamics of the random field Curie-Weiss model.} {\it Annals of Applied Probability} (2000), Vol. {\bf 10}, No 4, 1212 -- 1245. 

\bibitem{GDH} {\sc Greven, A.} and {den Hollander, F.} {\rm Large deviations for a random walk in random environment.}
{\it Ann. Probab.} (1994), Vol. {\bf 22}, No 3, 1381--1428.

\bibitem{guiliv} {\sc Guillotin-Plantard, N.} and {\sc Schott, R.} {\it Dynamic random walks: Theory and applications.} {\rm Elsevier} (2006).

\bibitem{gui1} {\sc Guillotin, N.} {\rm Asymptotics of a dynamic random
walk in a random scenery I. A law of large numbers.} {\it Annales de l'Institut Henri Poincar\'{e} - Probabilit\'{e}s et Statistiques} (2000), Vol. {\bf 36}, No 2, 127--151.

\bibitem{gui2} {\sc Guillotin, N.} {\rm Asymptotics of a dynamic random walk in a random scenery II. A functional limit theorem.} {\it Markov Processes and Related Fields} (1999), Vol. {\bf 5}, No 2, 201--218.

\bibitem{guiJTP} {\sc Guillotin-Plantard, N.} {\rm Dynamic $\bbZ^{d}$-random walks in a random scenery: a strong law of large numbers.} {\it J. Theoret. Probab.} (2001), Vol. {\bf 14}, No 1, 241--260.

\bibitem{arno} {\sc Guillotin-Plantard, N.} and {\sc Le Ny, A.} {\rm Transient random walks on 2d-oriented lattices.} {\it Theory of Probability and its Applications} (2006). To appear.

\bibitem{GPS} {\sc Guillotin-Plantard, N.} and {\sc Schott, R.} {\rm Distributed algorithms with dynamical random transitions.} {\it Random Structures and Algorithms} (2002), Vol. {\bf 21}, No 3-4, 376 -- 395.

\bibitem{KPZ} {\sc Koukiou, F.},  {\sc Pétritis, D.} and {\sc Zahradnik, M.} {\rm Extension of the Pirogov-Sinai theory to a class of quasiperiodic interactions.} {\it Commun. Math. Phys.} (1988), Vol. {\bf 118}, 365 -- 383. 

\bibitem{kui}{\sc Kuipers, L.} and {\sc Niederreiter, H.} {\em Uniform distribution of sequences.} {\rm Wiley and sons (1974).}

\bibitem{Kul} {\sc Külske, C.} {\rm Metastates in disordered mean field models: random field and Hopfield models.} {\it J. Statist. Phys.} (1997), Vol. {\bf 88}, 1257 -- 1293.

\bibitem{nie} {\sc Osgood F., C.} {\it Diophantine approximation and its applications.}
{\rm Academic Press (1973).}

\bibitem{Var} {\sc Varadhan, S.R.S.} {\rm Asymptotic probabilities and differential equations.} {\it Comm. on pure and applied math.}
 (1966), Vol. {\bf 19}, 261--286.

\bibitem{yo}  {\sc Yoccoz, J-C.} {\rm Sur la disparition de la propri\'{e}t\'{e} de Denjoy-Koksma en dimension 2.}
 {\it Ast\'{e}risque,} {\bf 231,} (1995).

\end{thebibliography}
\end{document}